\documentclass
[12pt]{amsart}

\RequirePackage[utf8]{inputenc}
\RequirePackage{xcolor}
\RequirePackage{hyperref}
\RequirePackage{amsfonts,amsmath,amsthm,amssymb}
\RequirePackage{mathtools} 
\RequirePackage{easybmat}  
\RequirePackage{comment}
\RequirePackage{multirow}  
\RequirePackage{xspace}
\RequirePackage{enumerate}
\RequirePackage{algorithm}
\RequirePackage{rotating}
\RequirePackage{algpseudocode}
\RequirePackage{nicematrix}
\RequirePackage{mathrsfs}
\RequirePackage{dynkin-diagrams}
\RequirePackage{rank-2-roots}
\RequirePackage{pgf,tikz,pgfplots,tikz-cd}
\usetikzlibrary{arrows,arrows.meta,decorations.markings,decorations.pathreplacing,calc}
\RequirePackage{tcolorbox}

\usepackage{ifdraft}
\usepackage{parskip}
\usepackage{tabularray}
\usepackage[obeyDraft]{todonotes}

\ifdraft{%
    \usepackage[
        margin=1.1in,
        left=0.2in,
        right=2in,
        marginparsep=0.2in,
        marginparwidth=1.7in,
    ]{geometry}
    \usepackage{soul}
    \sethlcolor{orange!20}
    \setuptodonotes{%
        size=\footnotesize,
        color=orange!20,
    }
}{
    \usepackage[margin=1.1in]{geometry}
    
}

\pgfplotsset{compat=1.15}

\DeclareMathOperator{\Aut}{Aut}

\DeclareMathOperator{\ad}{ad}

\DeclareMathOperator{\Ker}{Ker}
\DeclareMathOperator{\Id}{Id}
\newcommand{\ip}[1]{\langle #1 \rangle}

\numberwithin{equation}{section}

\newtheorem{thm}{Theorem}[section]
\newtheorem*{thm*}{Theorem}

\newtheorem*{prop*}{Proposition}

\newtheorem{question}[thm]{Question}
\newtheorem*{question*}{Question}

\newtheorem{cor}[thm]{Corollary}

\newtheorem{lemma}[thm]{Lemma}

\newtheorem*{conj*}{Conjecture}
\theoremstyle{definition}
\newtheorem{defin}[thm]{Definition}

\newtheorem{example}[thm]{Example}
\newtheorem{remark}[thm]{Remark}
\newtheorem*{remark*}{Remark}
\newtheorem*{remarks*}{Remarks}
\newtheorem*{notation*}{Notation}

\hypersetup{
    colorlinks=true,
    linkcolor=blue,
    filecolor=magenta,      
    urlcolor=cyan,
}

\newcommand{\iptnR}{\ip{ \ , \ }_{\mathfrak t(n,\mathbb R)}}

\begin{document}

\title[Einstein solvmanifolds as submanfiolds]{Einstein solvmanifolds as submanifolds of symmetric spaces}
\author[M. Jablonski]{Michael Jablonski}
\thanks{This work was supported in part by NSF grant DMS-1612357.}
\maketitle

\begin{abstract}  Lying at the intersection of  Ado's theorem and the Nash embedding theorem,  we consider the problem of finding faithful representations of Lie groups which are simultaneously isometric embeddings.  Such special maps are found for a certain class of solvable Lie groups which includes all Einstein and Ricci soliton solvmanifolds, as well as all Riemannian 2-step  nilpotent Lie groups.  As a consequence, we extend work of Tamaru by showing that all Einstein solvmanifolds  can be realized as submanifolds (in the submanifold geometry) of a symmetric space.
\end{abstract}

Recently, Tamaru \cite{Tamaru:ParabolicSubgroupsOfSemisimpleLieGroupsAndEinsteinSolvmanifolds} constructed a collection of explicit examples of Einstein solvmanifolds which are isometrically embedded into symmetric space.  There, for each symmetric space of non-compact type, a finite number of  examples are produced; these examples are also interesting in that they are not totally geodesic.   Naturally, it raises the question of whether or not all Einstein solvmanifolds arise this way.

\begin{question*}  Does every Einstein solvmanifold arise as a  submanifold of a symmetric space? 
\end{question*}

To approach this question, we reframe it as a question on solvable Lie groups.   Recall that every symmetric space is a totally geodesic submanifold of the  symmetric space $GL(n,\mathbb R)/ O(n,\mathbb R)$, which can be naturally identified with the group of lower triangular matrices $T(n,\mathbb R)\subset GL(n,\mathbb R)$ endowed with a  left-invariant, symmetric metric.  This symmetric space is a product of a 1-dimensional Euclidean de Rham factor and an irreducible, non-flat factor which is Einstein; the symmetric metric is unique up to scaling on this non-flat factor.  Even stronger than being a submanifold of a symmetric space, those examples by Tamaru can all be realized as subgroups of the group of lower triangular matrices.  As every Einstein solvmanifold is isometric to a completely solvable Lie group with left-invariant metric 
\cite{Lauret:EinsteinSolvmanifoldsAreStandard}, it seems natural to ask the even more delicate question of whether one can achieve an embedding of Einstein solvmanifolds  via faithful representations into the lower triangular matrices.

\begin{defin}  Consider the group of lower triangular matrices $T(n,\mathbb R)\subset GL(n,\mathbb R)$ endowed with a fixed choice of left-invariant, symmetric metric.  Let $S$ be a solvable Lie group with left-invariant metric.  By an isomorphic, isometric embedding, we mean a faithful representation
	$$\phi : S \to T(n,\mathbb R)$$
which is an isometric embedding of Riemannian manifolds.
\end{defin}

\begin{thm}\label{thm: main theorem for Einstein}  Let $S$ be a 
completely solvable Lie group endowed with a left-invariant Einstein metric.  Then there exists $n \in \mathbb N$ and an isomorphic, isometric embedding $\phi : S \to T(n,\mathbb R)$ where $T(n,\mathbb R)$ is endowed with a left-invariant, symmetric metric.
\end{thm}

Essentially, what we are describing above is a natural combination  of Ado's Theorem (which says every Lie group is locally linear) and Nash's Embedding Theorem for Riemannian manifolds.  On a solvable Lie group, if an Einstein metric exists it is unique (up to scaling and isometry) \cite{Heber,Lauret:EinsteinSolvmanifoldsAreStandard} 
and so the symmetric metric on $T(n,\mathbb R)$ plays the natural r\^ole of the Eucliean metric  on $\mathbb R^n$.  Furthermore, Einstein metrics are known to have maximal symmetry and  this is even more reason to consider them as a preferred choice of metric on $T(n,\mathbb R)\cap SL(n,\mathbb R)$ \cite{GordonJablonski:EinsteinSolvmanifoldsHaveMaximalSymmetry}.

One consequence of the theorem above is an affirmative answer to the question at the beginning.

\begin{cor}  Every 
Einstein solvmanifold can be realized as a submanifold of a symmetric space.
\end{cor}

Note, by a recent result of B\"ohm-Lafuente \cite{Bohm-Lafuente:NonCompactEinsteinManifoldsWithSymmetry}, every non-compact, homogeneous Einstein manifold is an Einstein solvmanifold.  In this way, we see that all non-compact, homogeneous Einstein manifolds arise as submanifolds of symmetric spaces.

It seems natural to ask if any Lie group  with left-invariant metric  can be realized as a submanifold  (in the submanifold geometry) of some canonical Riemannian Lie group.  In some special cases, we can say that this is true.

\begin{cor} \label{cor: main result for ricci solitons}
Let $S$ be a simply-connected, completely solvable Lie group endowed with a left-invariant Ricci soliton metric.  Then there exists $n \in \mathbb N$ and an isomorphic, isometric embedding $\phi : S \to T(n,\mathbb R)$ where $T(n,\mathbb R)$ is endowed with a left-invariant, symmetric metric.
\end{cor}

Recall, for solvable Lie groups that cannot admit an Einstein metric, such as nilpotent Lie groups, a Ricci soliton is arguably the preferred geometry  \cite{Lauret:StandardEinsteinSolvAsCriticalPoints,GordonJablonski:RicciSolitonSolvmanifoldsHaveInfinitesimalMaximalSymmetry}.  Unlike in the Einstein setting where solvmanifolds must be simply-connected, we must add the hypothesis of simple-connectivity in the Ricci soliton setting to insure that the isomorphism $\phi$ exists at the Lie group level and not just the Lie algebra level.

For a general metric on a solvable Lie group, we do not know if such embeddings exist.  However, we can say more for 2-step nilpotent Lie groups.

\begin{thm}  \label{thm:  main result for 2-step nilpotent}
Every simply-connected, 2-step nilpotent Lie group with left-invariant metric  isomorphically, isometrically embeds into some $T(n,\mathbb R)$ endowed with a left-invariant, symmetric metric.
\end{thm}

We suspect the same holds for  every nilpotent group with left-invariant metric, but do not have a proof at hand.

Finally, we comment on the compact setting to put our results above in perspective.  The natural candidate in which to embed a compact Lie group   is $O(n,\mathbb R)$ with its bi-invariant metric.  However, the above results cannot hold in this setting.  Any subgroup of $O(n,\mathbb R)$ would pickup a bi-invariant metric itself with the submanifold geometry, and there are many metrics on a given non-abelian, compact Lie group that are not bi-invariant.  

Perhaps there is another, natural candidate besides the bi-invariant metric that one could use for $O(n,\mathbb R)$; this would be interesting to know.

\subsubsection*{Acknowledgments}  The author would like to thank Megan Kerr and Tracy Payne for conversations in which the question above arose.

\section{Preliminaries}\label{sec: preliminaries}
Our first observation is that one can reduce to the study of metric Lie algebras, that is Lie algebras endowed with an inner product.  We explain this below.  Note, throughout the paper, we will denote Lie groups by their capital, Roman letter and the corresponding Lie algebras by their lower case, mathfraktur letter.

Recall, simply-connected Lie groups with left-invariant metrics are in one-to-one correspondence with Lie algebras with inner products. If one has an isomorphic, isometric embedding $\phi: S \to T(n,\mathbb R)$, then one immediately has the induced homomorphism between their Lie algebras
	$$\phi :  \mathfrak s \to \mathfrak t(n,\mathbb R)$$
which is an isomorphism onto its image and which preserves the inner product on the Lie algebra, i.e.
	$$\ip{X,Y} _ \mathfrak s = \ip{\phi_*(X), \phi_*(Y)   } _ {\mathfrak t(n,\mathbb R)}   \quad  \mbox{   for  } X,Y\in\mathfrak s.$$
Conversely, having such a map at the Lie algebra level lifts to a homomorphism between the simply-connected Lie groups $S$ and $T(n,\mathbb R)$ which is an isomorphism onto its image.  Being a homomorphism, we see that the left-invariant metric on $T(n,\mathbb R)$ pulls back to a left-invariant metric on $S$; that is, the map $\phi$ is an isomorphic, isometric embedding.

We will often simply denote either of the inner products on $\mathfrak s$ or $\mathfrak t(n,\mathbb R)$ by $\ip{ \ , \ }$ as long as the context is clear.  Further, we will denote by $\phi^*\ip{ \ , \ }$ the induced inner product on $\mathfrak s$ that one gets by restricting the inner product on $\mathfrak t(n,\mathbb R)$ to the image $\phi_*(\mathfrak s)$; at the group level this corresponds to the submanifold geometry of $\phi(S)$ sitting inside of $T(n,\mathbb R)$.

\subsection*{The metric on $T(n,\mathbb R)$}
The group $T(n,\mathbb R)$ has two seemingly natural left-invariant metrics which we describe via the corresponding inner products on $\mathfrak t(n,\mathbb R)$.  First is the Frobenius inner product.  This is just the restriction of the usual inner product on $\mathfrak{gl}(n,\mathbb R)$ given by
	\begin{equation}\label{eqn: fixed Einstein metric on T(n,R)}
	B(X,Y) = tr \ (XY^t)  \quad \mbox{ for } X,Y\in\mathfrak{gl}(n,\mathbb R).
	\end{equation}
Although familiar, there is no geometric justification for using this metric.

Using the Frobenius inner product, we can describe the symmetric inner product on $\mathfrak t(n,\mathbb R)$.  (We call it the symmetric inner product as the corresponding left-invariant metric on $T(n,\mathbb R)$ is symmetric.)

Write $\mathfrak t(n,\mathbb R) = \mathfrak a + \mathfrak n$ where $\mathfrak a$ is the set of diagonal matrices and $\mathfrak n$ is the set of strictly lower triangular matrices.  The symmetric inner product on $\mathfrak t(n,\mathbb R)$ which we will work with is  given by
	\begin{equation}\label{eqn: standard iptnR}
	\begin{alignedat}{2}
	\iptnR  &=   2 B \quad & \mbox{ on } \mathfrak a \times \mathfrak a \\
	\iptnR &= B \quad & \mbox{ on } \mathfrak n\times \mathfrak n
	\end{alignedat}
	\end{equation}
with $\mathfrak a$ and $\mathfrak n$ orthogonal.  The corresponding left-invariant geometry on $T(n,\mathbb R)$ is isometric to $GL(n,\mathbb R)/O(n,\mathbb R)$
 with the symmetric metric and so stands out geometrically.

Our main results hold for both geometries on $T(n,\mathbb R)$, but we are primarily interested in the symmetric space geometry (or equivalently, the Einstein geometry on $T(n,\mathbb R)\cap SL(n,\mathbb R)$).

\subsection*{Picking a scale}

Up to scaling the non-flat factor, there is only one symmetric metric on $T(n,\mathbb R)$.  However, for the above theorems, one actually has to pick a scale.  Instead, one can equivalently look for representations $\phi : \mathfrak s \to \mathfrak t(n,\mathbb R)$ which satisfy
	\begin{equation} \label{eqn: embedding up to scale}
	\phi^* \ip{ \ , \ }_{\mathfrak t(n,\mathbb R)} = c \ip{ \ , \ }_\mathfrak s   \quad \mbox{ for some } c\in\mathbb R,
	\end{equation}
where the inner product $\iptnR$ is the one given in Equation \ref{eqn: standard iptnR}.  
Then one simply chooses the ``right'' symmetric metric on the left so as to make $c=1$ on the right.  In the work that follows, we will find $\phi$ which satisfies Eqn.~\ref{eqn: embedding up to scale}.

\begin{remark}
It is hard to justify any one scale from another and so we do not.  If one wanted to fix a scale, then  one could fix the scalar curvature on $T(n,\mathbb R)\cap SL(n,\mathbb R)$ to be -1.  Or one could instead embed into $T(n+1,\mathbb R)\cap SL(n+1,\mathbb R)$ to eliminate the Euclidean de Rham factor all together and then fix the metric with a particular value of scalar curvature.
\end{remark}

\subsection*{Structure of Einstein solvmanifolds}

Let $S$ be a solvable Lie group which admits an Einstein metric with Lie algebra $\mathfrak s$.  The Lie group $S$ is necessarily simply-connected.  The Lie algebra decomposes as $\mathfrak s = \mathfrak a + \mathfrak n$, where $\mathfrak n$ is the nilradical and $\mathfrak a$ is an abelian subalgebra acting reductively on $\mathfrak n$. Further, under the Einstein metric, $\mathfrak a \perp \mathfrak n$.

From here going forward, we assume that $S$ is completely solvable - i.e., $\ad\, X$ has real eigenvalues for $X\in\mathfrak s$.  Up to isometry, this captures the full class of Einstein solvmanifolds.  For such an $\mathfrak s$ with Einstein metric, it is known that there exists some $A\in\mathfrak a$ such that $D=\ad\, A$ has positive eigenvalues on $\mathfrak n$; consequently, $\mathfrak s$ has no center.  Furthermore, $\ad\,\mathfrak a$ acts by symmetric transformations (relative to the Einstein inner product).

For a detailed discussion of the structure of Einstein solvmanifolds, we refer the interested reader to \cite{Lauret:EinsteinSolvmanifoldsAreStandard}.

\section{Isomorphically embedding Einstein solvmanifolds}

To prove our main results, we  work in a slightly more general setting.  

\begin{defin}\label{defin: admissible metric Lie algebra}  We say a metric Lie algebra $(\mathfrak s , \ip{ \ , \ })$ is admissible if it satisfies the following conditions:
	\begin{enumerate}
	\item  $\mathfrak s = \mathfrak a + \mathfrak n$  is completely solvable,
	\item  $\mathfrak n$ is the nilradical of $\mathfrak s$,
	\item $\mathfrak a$ is an abelian subalgebra, 
	\item for each $X\in\mathfrak a$, $\ad X:\mathfrak s\to\mathfrak s$ is symmetric relative to the given inner product,
	\item there is some $A\in\mathfrak a$ such that $D=\ad\,A$ has positive eigenvalues on $\mathfrak n$.
	\end{enumerate}
\end{defin}

\begin{remark}
For each $X\in\mathfrak a$, the map $\ad X:\mathfrak s\to\mathfrak s$ being  symmetric implies that it is diagonalizable over $\mathbb R$.  As $\mathfrak a$ is abelian, we have that the collection of maps $\ad \mathfrak a$ are simultaneously diagonalizable and their weght spaces are orthogonal.  In particular, condition (v)  implies the zero weight space of the $\ad\,\mathfrak a$-action on $\mathfrak s$ is $\mathfrak a$ itself and  so   $\mathfrak a \perp \mathfrak n$.
\end{remark}

\begin{thm} \label{thm: main general theorem}
Let $S$ be a simply-connected, completely solvable Lie group with left-invariant metric whose metric Lie algebra is admissible. 
There exists an isomorphic, isometric embedding $\phi: S \to T(n,\mathbb R)$, up to scaling, where $T(n,\mathbb R)$ is the group of lower triangular matrices endowed with the left-invariant, symmetric metric given in Eqn.~\ref{eqn: fixed Einstein metric on T(n,R)}.
\end{thm}

\begin{remark}  Given the structure of Einstein solvmanifolds outlined in Section \ref{sec: preliminaries}, clearly the above theorem implies Theorem \ref{thm: main theorem for Einstein}. The more general case of Ricci soliton solvmanifolds and general Riemannian, 2-step nilpotent Lie groups fall under this frame work, as well - see Section \ref{sec: beyond Einstein metrics}.
\end{remark}

The eigenspace decomposition of $D$ on $\mathfrak n$  gives a grading to the nilradical 
	$$\mathfrak n =  \mathfrak n_{\lambda_1} \oplus \dots \oplus \mathfrak n_{\lambda_k},$$
where $0<\lambda_1 < \dots < \lambda_k$ are the eigenvalues of $D$.  
By grading, we mean that the above satisfies $[\mathfrak n_{\lambda_i},\mathfrak n_{\lambda_j}] \subset \mathfrak n_{\lambda_i+\lambda_j}$.

\begin{remark}  Given the grading above, we observe that $\mathfrak n_{\lambda_k} \subset \mathfrak z(\mathfrak n)$, the center of $\mathfrak n$.  Additionally, if we consider the $r$-largest eigenvalues $\{\lambda_{k-r+1},\dots,\lambda_k\}$ of $D$, then the sum of the corresponding weight spaces $\mathfrak n_{\lambda_{k-r+1}}\oplus\dots\oplus\mathfrak n_{\lambda_{k}}$ is an ideal of $\mathfrak s$.
\end{remark}

\subsection{Strategy}  We begin with an outline of the proof of Theorem \ref{thm: main general theorem}.  Step 1 is to control the metric on the nilradical, obtaining a representation which keeps the eigenspaces of $D$ orthogonal and simply rescales the inner product on each eigenspace.  Step 2 is to modify the representation so as to obtain the same scaling on all eigenspaces of the nilradical simultaneously.  Step 3 is to modify the representation so that it scales the inner product on the complementary subalgebra $\mathfrak a$ by the same factor, while keeping $\mathfrak a$ orthogonal to $\mathfrak n$.

\subsection{Tactics.}
To carry out our strategy, we start with a familiar representation and make various modifications.  Let $\mathfrak s$ be a solvable Lie algebra as above and consider the adjoint representation 
\begin{equation}\label{eqn: fundamental representation}
	 \phi = \ad : \mathfrak s \to \mathfrak{gl} (\mathfrak s).
\end{equation}
As $D = \ad\,A$ is   non-singular on $\mathfrak n $ for some $A\in\mathfrak a$, by hypothesis, $\mathfrak s$ has no center and we have that $\phi$ is a faithful representation.  

Consider $D=\ad\,A$ where the eigenvalues are positive on $\mathfrak n$.  Now choose an ordered  basis of $\mathfrak s$ which is a union of bases of the weight spaces of $\ad\,\mathfrak a$ and ordered so that the lower eigenspaces of $D$ come first.  In this way, we see that $\phi(\mathfrak s) \subset \mathfrak t(n,\mathbb R)$.

Using such a representation as a seed, we will modify it by  composing  with automorphisms and also adding such representations together.  Consider the following examples.

\begin{example}\label{ex: conj postcomposed with rep}
Given a faithful representation $\phi: \mathfrak s \to \mathfrak{t}(n,\mathbb R)$ and $L\in T(n,\mathbb R)$, we have a new faithful representation $C_L\circ \phi: \mathfrak s \to \mathfrak{t}(n,\mathbb R)$, where $C_L$ is conjugation by $L$.
\end{example}

\begin{example}\label{ex: precomp rep with aut}
Given a faithful representation $\phi: \mathfrak s \to \mathfrak{t}(n,\mathbb R)$ and an automorphism $a\in \Aut(\mathfrak s)$,  we have a new faithful representation $\phi\circ a: \mathfrak s \to \mathfrak{t}(n,\mathbb R)$.
\end{example}

\begin{example}\label{ex: sum of reps}
Given two representations $\phi_1 : \mathfrak s \to \mathfrak{gl}(V_1)$ and $\phi_2  : \mathfrak s \to \mathfrak{gl}(V_2)$, we may build a new representation 
	\begin{equation} \label{eqn: adding representations}
	\phi_1 + \phi_2 :  \mathfrak s \to \mathfrak{gl}(V_1)\oplus \mathfrak{gl}(V_2) \subset \mathfrak{gl}(V_1\oplus V_2).
	\end{equation}
If both representations are lower triangular, so is their sum. If either $\phi_i$ is faithful, so is $\phi_1+\phi_2$.

\end{example}

\begin{remark}\label{remark: coordinates for sums of reps}
As the coordinates of $\phi_1$ and $\phi_2$ lie in different blocks, they are orthogonal relative to   the Frobenius inner product.  Likewise, for representations into $\mathfrak t(n,\mathbb R) \times \mathfrak t(m,\mathbb R) \subset \mathfrak t(n+m,\mathbb R)$, the coordinates of $\phi_1$ and $\phi_2$ are orthogonal relative to the symmetric inner product.  Moreover, $| \phi_1(X) + \phi_2(X)    |^2   =  | \phi_1(X)   |^2  +  |  \phi_2(X)  |^2$.  These properties will be essential in the work that follows.
\end{remark}

\section{Step 1: Initial pass on the nilradical}

\begin{lemma} \label{lemma:step 1}
Let $(\mathfrak s, \ip{\ , \ }_\mathfrak s) $ be an admissible metric Lie algebra.  There exists a representation $\phi: \mathfrak s \to\mathfrak t (n,\mathbb R)$ such that the inner product $\phi^* \iptnR$ on $\mathfrak s$ satisfies 
	\begin{enumerate}
	\item the eigenspaces of $D$ are orthogonal and
	\item on each eigenspace of $D$ on $\mathfrak n$, the inner product $\phi^* \iptnR$ is simply a rescaling of $\ip{ \ , \ }_{\mathfrak s}$.
	\end{enumerate}
\end{lemma}

\begin{remark} The inner product obtained will likely have different scalings on each eigenspace.  This is dealt with in Step 2.
\end{remark}

We now prove this lemma.  Consider the solvable algebra $\mathfrak s =  \mathfrak a \ltimes \mathfrak n$.  As $\mathfrak s$ has no center, its adjoint representation 
	$$\phi  = \ad   : \mathfrak s \to \mathfrak{gl}(\mathfrak s).$$
is faithful.

Recall, as a vector space $\mathfrak s$ decomposes as a direct sum of the eigenspaces of $D$
	\[\mathfrak s = \mathfrak a + \mathfrak n_{\lambda_1} + \dots + \mathfrak n_{\lambda_k},\]
where $\mathfrak a = \Ker D$ and $0<\lambda_1< \dots <\lambda_k$.  For the sake of convenient notation going forward, we will write $\mathfrak a = \mathfrak n_0$ as $\mathfrak a$ is the zero eigenspace.  Warning, this notation is slightly misleading as $\mathfrak a$ is not contained in the nilradical $\mathfrak n$.

Next, we will choose a particular basis of $\mathfrak s$ so that $\phi(\mathfrak s)$ is in the set of lower triangular matrices $\mathfrak t(n,\mathbb R)$.  Denote the weights of the $\ad\mathfrak a$-action on $\mathfrak s$ by $\{\alpha\}$; to each weight $\alpha$ we have the corresponding weight space $\mathfrak n_\alpha$.  Note, the zero weight space of the action is still $\mathfrak a$ since $D$ has all positive eigenvalues on $\mathfrak n$.

As $\ad\mathfrak a$ preserves each eigenspace $\mathfrak n_{\lambda_i}$ of $D$, we have $\mathfrak n_{\lambda_i}=\bigoplus_\alpha \left(\mathfrak n_{\lambda_i}\cap \mathfrak n_\alpha   \right)$.  We choose an (ordered) orthonormal basis  $\{e_1,\dots, e_n\}$ of $\mathfrak s$    which consists of a union of bases of the subspaces $\mathfrak n_{\lambda_i}\cap \mathfrak n_\alpha$, which are ordered to have the lesser  eigenvalues of $D$ first.  (We impose no ordering relative to the other weights.)  In this way, we have $\phi(\mathfrak s)\subset \mathfrak t(n,\mathbb R)$.  Notice that $\phi(\mathfrak a)$ consists of diagonal matrices while $\phi(\mathfrak n)$ consists of strictly lower triangular matrices.

For future use, we go ahead analyze $\phi(\mathfrak n_\alpha)$.  Recall, using the inner product on $\mathfrak s$, we have the Frobenius inner product on $\mathfrak{gl}(\mathfrak s)$ given by
	\begin{equation}\label{eqn: frobenius norm using basis}
	\ip{A,B} = tr (AB^t) = \sum_i \ip{Ae_i,Be_i}_\mathfrak s,
	\end{equation}
for $A,B\in\mathfrak{gl}(\mathfrak s)$ and an orthonormal basis $\{e_i\}$ of $\mathfrak s$.  We consider the corresponding symmetric metric on $T(n,\mathbb R)$ as above, see Section \ref{sec: preliminaries}.

Given two weight spaces $\mathfrak n_\alpha$ and $\mathfrak n_\beta$, we see that $[\mathfrak n_\alpha,\mathfrak n_\beta] \subset\mathfrak n_{\alpha+\beta}$.  Now assume $\alpha\not=\beta$ and let $\gamma$ be another weight of the $\ad\mathfrak a$-action on $\mathfrak s$.  Take $X_\alpha\in\mathfrak n_\alpha$ and  $X_\beta\in\mathfrak n_\beta$.  Since $\phi(X_\alpha): \mathfrak n_\gamma \to \mathfrak n_{\alpha + \gamma}$, $\phi(X_\beta): \mathfrak n_\gamma \to \mathfrak n_{\beta + \gamma}$, and $\alpha+\gamma \not= \beta + \gamma$, we see that the maps $\phi(X_\alpha)$ and $\phi(X_\beta)$ share no common coordinates.  Computing the inner product using Eqn.~\ref{eqn: frobenius norm using basis}, we see that 
$$\phi (\mathfrak n_\alpha) \perp \phi(\mathfrak n_\beta),$$
with respect to the symmetric inner product on $\mathfrak t(n,\mathbb R)$.

\subsection{The inner product on the highest eigenspace $\mathfrak n_{\lambda_k}$}
Consider $\mathfrak n_{\lambda_k}$, the eigenspace of $D$ where $\lambda_k$ is the largest eigenvalue of $D$.  As $\ad \mathfrak a$ commutes with $D$, we have $\mathfrak n_{\lambda_k}~=~\bigoplus_{\alpha}~\left(\mathfrak n_{\lambda_k}\cap\mathfrak n_\alpha\right)$, where $\{\alpha\}$ are the weights of the $\ad\mathfrak a$-action. 
Take  $e_i \in \mathfrak n_{\lambda_k}\cap\mathfrak n_\alpha$.  Since $e_i$ is central in $\mathfrak n$, we are able to easily compute $\phi(e_i)$, as follows.

Let $a = \dim \mathfrak a$ and $\{A_1,\dots, A_a\}$ be our orthonormal basis of $\mathfrak a$.  Then 
	\[\phi(e_i) =  \sum_{1\leq j \leq a } -\alpha (A_j)  E_{i , j}\]
Observe, the adjoint representation is then a  faithful  representation with the property that  it rescales the metric on the weight spaces of $\ad\,\mathfrak a$ acting on $\mathfrak n_{\lambda_k}$.  More precisely, on each weight space we have
	\begin{equation}\label{eqn: inner product on image on weighspaces of center}
	\phi^* \iptnR = c_\alpha \ip{ \ , \ }_\mathfrak s,
	\end{equation}
where $c_\alpha = { \sum \alpha (A_j)^2}$.  Note, $c_\alpha$ cannot be zero since $D$ has no kernel on $\mathfrak n_{\lambda_k}$.

Next, we adjust our representation so that these weight spaces rescale by the same constant. Consider a   matrix $L\in T(n,\mathbb R)$   which is the identity on $\mathfrak a \oplus \mathfrak n_{\lambda_1} \oplus \dots \oplus \mathfrak n_{\lambda_{k-1}}$, while on each weight space $\mathfrak n_{\lambda_k}\cap\mathfrak n_\alpha$ is a lower triangular matrix $L_\alpha$.

Let $C_L$ denote conjugation by $L$ in $\mathfrak{gl}(\mathfrak s)$.  The composition  $C_L \circ \phi$ is still a faithful representation, cf.\ Example \ref{ex: conj postcomposed with rep}.  The images of the weight spaces  $\mathfrak n_{\alpha}$ of $\ad\,\mathfrak a$ are still orthogonal relative to the symmetric inner product, as $L$ preserves each weight space.

Writing out block matrices, it is quick to see that, for $e_i\in \mathfrak n_{\lambda_k}\cap\mathfrak n_\alpha$, we have $\phi(e_i)L^{-1}~=~\phi(e_i)$.  This yields
	\[(C_L\circ \phi )(e_i)
		=L\phi(e_i) =   \sum_{1\leq j \leq a } -\alpha (A_j)  L_\alpha E_{i , j}.\]
Here, we are naturally identifying each column vector $E_{i,j}$ with the vector $e_i\in\mathfrak n_{\lambda_k}\cap\mathfrak n_\alpha\subset\mathfrak n$.

Taking $e_{i_1},e_{i_2}\in\mathfrak n_{\lambda_k}\cap\mathfrak n_\alpha$, we compute
	\begin{align*}
	(C_L\circ \phi)^*\ip{e_{i_1},e_{i_2}}_{\mathfrak t(n,\mathbb R)} &= \ip{(C_L\circ \phi) e_{i_1} ,(C_L\circ \phi) e_{i_2}} _{\mathfrak t(n,\mathbb R)}\\
		&= \ip{\sum_{1\leq j_1 \leq a } -\alpha (A_{j_1})  L_\alpha E_{i_1 , j_1} , \sum_{1\leq j_2 \leq a } -\alpha (A_{j_2})  L_\alpha E_{i_2 , j_2} }_{\mathfrak t(n,\mathbb R)}\\
		&=\sum_{1\leq j \leq a } \alpha (A_{j})^2 \ip{  L_\alpha E_{i_1 , j},L_\alpha E_{i_2 , j} }_{\mathfrak t(n,\mathbb R)}\\
		&=\sum_{1\leq j \leq a } \alpha (A_{j})^2 \ip{  L_\alpha e_{i_1},L_\alpha e_{i_2} }_{\mathfrak s}
	\end{align*}

This shows that on the weight space $\mathfrak n_{\lambda_k}\cap\mathfrak n_\alpha$, we  have
	\begin{equation}\label{eqn: conjugating by L on the last eigenspace of D}
	(C_L\circ \phi)^*\iptnR = c_\alpha  \ip{L_\alpha\ ,L_\alpha \ }_\mathfrak s.
	\end{equation}

Take any $c>0$.  Choosing $L_\alpha = \sqrt{c/c_\alpha} \Id$ we see that we have a representation which satisfies 
	\begin{equation}\label{eqn: rep that scales on all of last eigenspace}
	(C_L\circ \phi)^*\iptnR = c \ip{\ , \ }_\mathfrak s
	\end{equation}
on all of $\mathfrak n_{\lambda_k}$.

Next, we adjust the representation further to take care of the remaining eigenspaces  $\mathfrak n_{\lambda_1} \oplus \dots \oplus \mathfrak n _{\lambda_{k-1}}$.

\subsection{The inner product on the remaining eigenspaces}

For a vector space $V$ with inner product $\ip{ \ , \ }$, recall that every symmetric, positive definite, bilinear form on $V$ can be represented  as 
	$$\ip{L(\cdot ),L(\cdot )}  \quad \mbox{ for } L\in GL(V) \mbox{ a lower triangular matrix. } $$ 

\begin{remark}\label{remark: adjusting ip on last eigenspace using lower triangular}
In light of Eqn. \ref{eqn: conjugating by L on the last eigenspace of D},  we can alter our faithful representation $\phi$ to have any prescribed inner product on $\mathfrak n_{\lambda_k}$, so long as the weight spaces of the $\ad\,\mathfrak a$ action remain orthogonal.
\end{remark}
We  apply this observation below to a related collection of admissible Lie algebras.

Consider the orthogonal eigenspace decomposition of $\mathfrak s$ relative to the derivation $D$
	\[\mathfrak s =\mathfrak a \oplus \mathfrak n_{\lambda_1} \oplus \dots \oplus \mathfrak n_{\lambda_k}.\]
As $\mathfrak n_{\lambda_k}$ is an ideal, we may consider the quotient algebra $\mathfrak s^{(2)} = \mathfrak s / \mathfrak n_{\lambda_k} \simeq \mathfrak a \oplus \mathfrak n_{\lambda_1} \oplus \dots \oplus \mathfrak n_{\lambda_{k-1}}$.  We endow $\mathfrak s^{(2)}$ with the inner product from $\mathfrak a \oplus \mathfrak n_1 \oplus \dots \oplus \mathfrak n_{k-1}$.  This metric Lie algebra is rarely Einstein, but it is admissible.  Note, the nilradical of $\mathfrak s^{(2)}$ is $\mathfrak n^{(2)} \simeq \mathfrak n_1 \oplus \dots \oplus \mathfrak n_{k-1}$.

Now consider the adjoint representation for $\mathfrak s^{(2)}$
	\[\phi^{(2)}:\mathfrak s^{(2)}\to \mathfrak{gl}(\mathfrak s^{(2)})   .\]
This is faithful as the induced action of $D$ on $\mathfrak n^{(2)}$ continues to have positive eigenvalues.  (Composing with the quotient of $\mathfrak s$ by $\mathfrak n_{\lambda_k}$, we may interpret $\phi^{(2)}$ as a representation of $\mathfrak s$.) As above, we choose an orthonormal basis of weight vectors for the $\ad \mathfrak a$ action on $\mathfrak n^{(2)}$.

For $c\in\mathbb R$, consider the symmetric, blinear form 
	$$  c\ip{ \cdot , \cdot  }_{\mathfrak s^{(2)}} - \ip{ \phi^{(2)}  \ \cdot , \phi^{(2)} \ \cdot  }_{\mathfrak t (n_2,\mathbb R)}     $$
on $\mathfrak a\oplus \mathfrak n_{\lambda_1}\oplus \dots \oplus \mathfrak n_{\lambda_{k-1}}$.  Here $n_2=\dim {\mathfrak s^{(2)}}$.  For $c\in\mathbb R$ large, this symmetric, blinear form is positive definite and with the property that  the weight spaces $\mathfrak n_\alpha$ of $\ad\,\mathfrak a$ are orthogonal.  
Applying Remark \ref{remark: adjusting ip on last eigenspace using lower triangular},  there is  a lower triangular matrix $L$ of $GL(\mathfrak a\oplus  \mathfrak n_{\lambda_1} \oplus \dots \oplus   \mathfrak n_{\lambda_{k-1}} )$ acting only on the last eigenspace   $\mathfrak n_{\lambda_{k-1}}$ 
such that
	\[  \ip{L \ \cdot ,L \ \cdot }_{\mathfrak s^{(2)}}  =    c\ip{ \cdot , \cdot  }_{\mathfrak s^{(2)}} - \ip{ \phi^{(2)}  \ \cdot , \phi^{(2)} \ \cdot  }_{\mathfrak{t}(n_2,\mathbb R)}\]
when restricted to the last eigenspace $\mathfrak n _{\lambda_{k-1}}$.

Let $\phi_1$ be the faithful representation given above (Eqn.~\ref{eqn: rep that scales on all of last eigenspace}) which simply rescales the metric on the highest eigenspace $\mathfrak n_{\lambda_k}$. 
Then the sum of representations 
	$$\phi_1 + C_L\circ \phi^{(2)}$$
has the property that it simply rescales on the last two eigenspaces $\mathfrak n_{\lambda_{k-1}}$ and  $\mathfrak n_{\lambda_k} $.  Note, one cannot expect that the rescalings on these two spaces are the same.  We take care of this in Step 2 of the proof of our main theorem.

\begin{remark} Although $\phi^{(2)}$ is not a faithful representation of $\mathfrak s$, since $\phi_1$ is faithful we have that $\phi = \phi_1  + C_L\circ \phi^{(2)}$ is faithful - cf.\ Example \ref{ex: sum of reps}.
\end{remark}

One finishes the proof of the lemma inductively.  Similar to the above, we may quotient out by the following ideals and endow the quotients with the corresponding inner products.  
	\begin{equation}\label{eqn: quotient algebras}
	\begin{split}
	\mathfrak s^{(3)} &= \mathfrak s /( \mathfrak n_{\lambda_{k-1}}\oplus \mathfrak n_{\lambda_k}) \simeq \mathfrak a \oplus \mathfrak n_{\lambda_1} \oplus \dots \oplus \mathfrak n_{\lambda_{k-2}}\\
	\mathfrak s^{(4)} &= \mathfrak s / ( \mathfrak n_{\lambda_{k-2}}\oplus \dots \oplus \mathfrak n_{\lambda_k})\simeq \mathfrak a \oplus \mathfrak n_{\lambda_1} \oplus \dots \oplus \mathfrak n_{\lambda_{k-3}}\\
	 & \ \ \vdots  \\
	\mathfrak s^{(k)} &= \mathfrak s / ( \mathfrak n_{\lambda_{2}}\oplus \dots \oplus \mathfrak n_{\lambda_k}) \simeq \mathfrak a \oplus \mathfrak n_{\lambda_1} 
	\end{split}
	\end{equation}
These metric algebras are all admissible.

Repeating the process above, we set $\phi = \phi_1 + C_L\circ \phi^{(2)}$.  Then we find the appropriate choice $L'$ such that 
	\[\phi + C_{L'}\circ \phi^{(3)}   \]
rescales the metrics on the last three eigenspaces.  Again, $\phi^{(3)}$ may be interpreted as a representation on $\mathfrak s$ by first taking the quotient $\mathfrak s \to \mathfrak s /( \mathfrak n_{\lambda_{k-1}}\oplus \mathfrak n_{\lambda_k}) = \mathfrak s^{(3)}$.

Continuing this process  for the remaining $\mathfrak s ^{(i)}$, we are able to build a representation such that the pullback metric simply rescales on each eigenspaces of $D$ on $\mathfrak n$, as desired.  This completes the proof of Lemma \ref{lemma:step 1}.

\section{Step 2:  Aligning all the scalings on the nilradical}

\begin{lemma}\label{lemma: step 2}
Let $(\mathfrak s, \ip{ \ , \ }_\mathfrak s)$ be an  admissible metric Lie algebra.  For each large $t\in\mathbb R$, there exists a representation $\phi_t: \mathfrak s \to\mathfrak t (n,\mathbb R)$ such that the inner product $\phi_t^* \iptnR$ on $\mathfrak s$ satisfies 
	\begin{enumerate}
	\item $\mathfrak a\perp \mathfrak n$, 
	\item restricted to $\mathfrak a \times \mathfrak a$ is constant (as a function of $t$), and 
	\item restricted to $\mathfrak n\times\mathfrak n$ equals 
	$t \ip{\ ,\ }_\mathfrak s$.
	\end{enumerate}
\end{lemma}

As in the previous section (see Eqn.\ref{eqn: quotient algebras}), we consider the quotient algebras with their natural metrics.  These spaces are all admissible metric Lie algebras.
	\begin{equation*}
	\begin{split}
	\mathfrak s^{(2)} &= \mathfrak s /\mathfrak n_{\lambda_k}\simeq \mathfrak a \oplus \mathfrak n_{\lambda_1} \oplus \dots \oplus \mathfrak n_{\lambda_{k-1}}\\
	\mathfrak s^{(3)} &= \mathfrak s / ( \mathfrak n_{\lambda_{k-1}}\oplus \dots \oplus \mathfrak n_{\lambda_k})\simeq \mathfrak a \oplus \mathfrak n_{\lambda_1} \oplus \dots \oplus \mathfrak n_{\lambda_{k-2}}\\
	 & \ \ \vdots  \\
	\mathfrak s^{(k)} &= \mathfrak s / ( \mathfrak n_{\lambda_{2}}\oplus \dots \oplus \mathfrak n_{\lambda_k}) \simeq \mathfrak a \oplus \mathfrak n_{\lambda_1} 
	\end{split}
	\end{equation*}
For notational convenience, we define $\mathfrak s^{(1)} :=\mathfrak s$.

Observe that $D$ induces a positive derivation on each $\mathfrak n^{(i)}$ whose eigenvalues are the first $k+1-i$ eigenvalues of $D$ on $\mathfrak n$.  For simplicity, we will denote each of these derivations by $D$.  Further, $a(t)=exp(tD)$ is an automorphism of $\mathfrak n^{(i)}$ which acts as scalar multiplication by $e^{t\lambda}$ on the eigenspace $\mathfrak n_{\lambda}$.

By Lemma \ref{lemma:step 1}, we may assume the existence of a faithful representation $\phi_i$ for each $\mathfrak s^{(i)}$  which keeps the eigenspaces of $D$ orthogonal and rescales the metric on each eigenspace $\mathfrak n_{\lambda}$ in the nilradical $\mathfrak n^{(i)}$.   As in Example \ref{ex: precomp rep with aut}, we may precompose with $a(t)$ and consider the representation $\phi_i\circ a(t)$.  For large $t$, we see that the scaling on the last eigenspace  $\mathfrak n_{\lambda_{k-i+1}}$ is the dominant scaling. Relatively speaking, the lower eigenspaces have their inner product essential shrunk as small as we like.

Now we consider 
	$$\phi = \phi_1 \circ a(t_1) + \dots + \phi_{k}\circ a(t_{k}),$$
similar to Example \ref{ex: sum of reps}.  As we will show, one may choose large $t_i$ so as to make all the eigenspaces scaled by the same amount.  

Fix $\ell\in\{1,\dots, k\}$.  For the inner product $\left(\phi_\ell\circ a(t_\ell) \right)^*\ip{\ , \ }$, denote the scaling constant for the eigenspace $\mathfrak n_{\lambda_i}$ by $c_{(\ell,{\lambda_i})}$.  Note, $c_{(\ell,{\lambda_i})}$ is a function of $t_\ell$.  As the eigenvalues of $D$ are ordered $\lambda_1<\lambda_2<\dots<\lambda_k$, we may increase $t_\ell$ until we have
	\[c_{(\ell,\lambda_1)}<c_{(\ell,\lambda_2)}<\dots<c_{(\ell,\lambda_\ell)}.\]
Further, 
this condition persists as $t_\ell$ grows, so it holds for large $t_\ell$.  Observe, for $i\leq \ell$, since $e^{t_\ell\lambda_i}/e^{t_\ell \lambda_{i-1}} \to \infty$ as $t_\ell\to\infty$, we have $c_{(\ell,\lambda_i)}-c_{(\ell,\lambda_{i-1})}\to\infty$.

For the inner product $\phi^*\ip{\ , \ }$, denote the scaling constant for the eigenspace $\mathfrak n_{\lambda_i}$ by $c_{\lambda_i}$.  Note, $c_{\lambda_i}$ is a function of $(t_1,\dots,t_k)$.  Employing Remark \ref{remark: coordinates for sums of reps}, we have 
	\[c_{\lambda_i} = c_{(1,\lambda_i)} + \dots + c_{(k-i+1,\lambda_i)}.\]
Increasing $t_2$, we may assume $c_{\lambda_1} < c_{\lambda_2}$.  This is because $c_{(2,\lambda_2)}-c_{(2,\lambda_1)}\to\infty$ as $t_2\to\infty$.  Likewise, increasing $t_3$, we may further assume $c_{\lambda_1} < c_{\lambda_2}<c_{\lambda_3}$.  Continuing through $t_k$, we maybe assume that
	\[  c_{\lambda_1} < c_{\lambda_2}<\dots < c_{\lambda_k} .\]
Observe, if one were to increase $t_\ell$, then one would still maintain the partial inequality of
	\[c_{\lambda_1} < c_{\lambda_2}<\dots<c_{\lambda_\ell}.\]
Finally, one can now increase $t_{k-1}$ until we have $c_{\lambda_{k-1}} =  c_{\lambda_k}$; that is, we have
	\[c_{\lambda_1} < c_{\lambda_2}<\dots<c_{\lambda_{k-1}} = c_{\lambda_k} .\]
By successively increasing $t_{k-2}, \dots, t_1$, we obtain
	\[c_{\lambda_1} = c_{\lambda_2} = \dots = c_{\lambda_k}.\]
That is, we have achieved $ \phi^*\iptnR =c \ip{ \ , \ }_\mathfrak s,$ 
for some $c>0$.  Going further, the above technique shows that for any $t\geq c$,  there is some choice of $(t_1,\dots,t_k)$ which solves  $ \phi^*\iptnR =t \ip{ \ , \ }_\mathfrak s$. This proves Lemma \ref{lemma: step 2}.

\section{Step 3:  adjusting the metric on $\mathfrak a$}

In this  section, we complete the proof of Theorem \ref{thm: main general theorem}.

Consider the vector space $W=\mathfrak a$.  Fixing an orthonormal basis $\{A_1,\dots,A_a\}$ of $\mathfrak a$, we may identify $\mathfrak a$ with the diagonal matrices $diag(W)$ in $\mathfrak{gl}(W)$.  In this way, the inner product on $\mathfrak a$ is precisely half the inner product on the diagonal using the symmetric inner product, cf.\ Eqns.~\ref{eqn: standard iptnR} \& \ref{eqn: frobenius norm using basis}.

Given  a linear map $L:\mathfrak a\to\mathfrak a$, we may view $L$ as a map $L:\mathfrak a \to diag(W)\subset\mathfrak{gl}(W)$.  As $diag(W)$ and $\mathfrak a$ are both abelian, the map $L$ is trivially an $\mathfrak a$-representation.  As the nilradical $\mathfrak n$ of $\mathfrak s$ is an ideal, we may view $L$ as a representation of $\mathfrak s$ by precomposing with the quotient $\mathfrak s \to \mathfrak s/\mathfrak n\simeq \mathfrak a$.  We  also denote this  representation of $\mathfrak s$ by $L:\mathfrak s\to\mathfrak{gl}(W)$.  Shortly, we will fix a useful choice of $L$.

Taking the representation $\phi_t : \mathfrak s \to \mathfrak{gl}(V)$ produced by Lemma \ref{lemma: step 2}, we form the sum of representations as in  Example \ref{ex: sum of reps}
	\[ \Phi_t = \phi_t + L: \mathfrak s \to \mathfrak{gl}(V) \oplus  \mathfrak{gl} (W)\subset \mathfrak{gl}(V\oplus W).   \]
With the correct choice of $L$, the representation $\Phi_t$ will be the sought after map satisfying the requirements in Theorem \ref{thm: main general theorem}.

\subsection*{Choosing $L:\mathfrak a \to \mathfrak a$}
On $\mathfrak a$, consider the symmetric, bilinear form $
\ip{ \phi_t \cdot , \phi_t \cdot   }$.  By construction of $\phi_t$, this bilinear form is constant in $t$.  For large $c\in\mathbb R$, we have
	$$ c\ip{ \cdot , \cdot } - \ip{ \phi_t \cdot , \phi_t \cdot   }$$
is a positive definite, symmetric, bilinear form  on $\mathfrak a$ - i.e., an inner product.  As such, there is an element $L\in GL(\mathfrak a)$ such that 
	$$\ip{ L \cdot , L\cdot }  =   c\ip{ \cdot , \cdot } - \ip{ \phi_t \cdot , \phi_t \cdot   }.$$
Note, $L=L(c)$ is a function of $c$ but not of $t$.  Using our identification of $\mathfrak a$ with $diag (W)\subset \mathfrak{gl}(W)$, we consider the map
	$$\Phi_t(A) = \phi_t(A) + L(A) \in \mathfrak{gl}(V) \oplus \mathfrak{gl}(W).$$
By construction, we have that for  large $c\in\mathbb R$ we can always achieve 
	\begin{equation}\label{eqn: solving the scaling on a x a}
		 \Phi_t^* \iptnR  =  c \ip{  \ , \  }_\mathfrak s \mbox{ on } \mathfrak a \times \mathfrak a \mbox{ for all large } t\in\mathbb R.
\end{equation}
Note, here $c$ and $t$ are independent of each other.

Choosing $c$ and $t$ large enough so that $\Phi_t$ exists and solves Eqn.~\ref{eqn: solving the scaling on a x a}, and taking  the max of these two, we may assume $c=t$ and now we have
	$$\Phi_t ^* \ip{  \ , \  }   = t \ip{  \ , \  }$$
on all of $\mathfrak s \times \mathfrak s$.   That is, we have produced an isomorphic, isometric embedding.  This proves the theorem.

\section{Beyond Einstein metrics}\label{sec: beyond Einstein metrics}
In the above, we prove the existence of an isomorphic, isometric embedding for a special class of  metrics, see Theorem \ref{thm: main general theorem}.  Although Ricci solitons do not fall into this class,  Corollary \ref{cor: main result for ricci solitons} follows immediately from the Einstein case as Ricci solitons can be realized as a subgroup with the submanifold geometry of an Einstein space \cite{Lauret:SolSolitons}.

To prove Theorem \ref{thm:  main result for 2-step nilpotent}, one simply needs to observe the well-known fact that every metric, 2-step nilpotent Lie algebra has a ``good'' rank 1 extension.  More precisely, let $\mathfrak n$ be 2-step nilpotent with center $\mathfrak z$.  Let $\mathfrak v$ be the orthogonal complement to $\mathfrak z$.  Then we may define a derivation $D$ of $\mathfrak n$ which is the identity on $\mathfrak v$ and twice the identity on $\mathfrak z$.  Observe, we can extend the inner product on $\mathfrak n$ to one on $\mathfrak s = \mathbb R D \ltimes \mathfrak n$ so that $D \perp \mathfrak n$.  In this way, $D$ is symmetric and the Lie algebra $\mathfrak s$ with the given inner product is admissible.  Denoting the corresponding simply-connected Lie groups by $N$ and $S$, we have that $N$ sits isometrically in $S$.  The claim follows from Theorem \ref{thm: main general theorem} applied to $S$.

\section{Final thoughts}

\begin{question}
In the proof of Theorem \ref{thm: main general theorem}, we showed that one can find $\phi$ which solves 
	$$\phi^*\ip{ \ , \ } = c \ip{ \ , \ },$$
for any large  $c\in\mathbb R$.  But can this be achieved for any $c\in\mathbb R$?
\end{question}

It is unclear what happens even in the Heisenberg case.  Using na\"ive attempts, one can easily get $c=1$, but it is unclear if a small value of $c$ is possible or if there is even a minimum.  It would be interesting to know if such a minimum has geometric significance.

Our construction above is likely far from efficient.

\begin{question}
What is the optimal dimension for achieving the above isomorphic, isometric embeddings?
\end{question}

Regarding the dimension of our target space, we note that if we fix $\dim S$ then there is a common $n = n(\dim S) \in\mathbb N$ such that an isomorphic embedding of $S$ into $T(n,\mathbb R)$ exists (up to scaling) for all $S$ with the given fixed dimension and satisfying the conditions of Theorem \ref{thm: main general theorem}.  This follows from noting that there are a finite number of partitions $\{\lambda_1, \dots, \lambda_k\}$ of $\dim S$ and that $n$ was determined by nothing more than the dimensions of $\mathfrak a$ and of the weight spaces $\mathfrak n_{\lambda_i}$ of the $\mathfrak a$-action.

\begin{remark}
What is not clear is whether or not there is a common value of $n$ after we fix a scale on the symmetric metric on $T(n,\mathbb R)$.  It would be interesting to know if this is the case.
\end{remark}

As a final remark, we note that in \cite{Kerr:NewExamplesOfNonsymmetricEinsteinSolv} Kerr constructs examples of Einstein solvmanifolds.  There it is claimed that those examples cannot arise as submanifolds of symmetric spaces, but as Kerr pointed out to us, what is proven is that those examples cannot arise as totally geodesic submanifolds, and so the results there are not at odds with the present work.

\providecommand{\bysame}{\leavevmode\hbox to3em{\hrulefill}\thinspace}
\providecommand{\MR}{\relax\ifhmode\unskip\space\fi MR }
\providecommand{\MRhref}[2]{%
  \href{http://www.ams.org/mathscinet-getitem?mr=#1}{#2}
}
\providecommand{\href}[2]{#2}

\end{document}